\documentclass{gtart}

\def\ifplaintex{\expandafter\ifx\csname documentclass\endcsname\relax}

\def\gtp{{\mathsurround=0pt\it $\cal G\mskip-2mu$eometry \&\ 
$\cal T\!\!$opology $\cal P\!$ublications}}  

\def\recd{{\small Received:\qua\receiveddate\ifx\reviseddate\relax
\else\qquad Revised:\qua\reviseddate\fi\par}} 

\def\lognumber#1{\def\thelognumber{#1}}
\def\volumenumber#1{\def\thevolumenumber{#1}}
\def\volumeyear#1{\def\thevolumeyear{#1}}
\def\papernumber#1{\def\thepapernumber{#1}}
\def\pagenumbers#1#2{\def\startpage{#1}\def\finishpage{#2}}
\def\published#1{\def\publishdate{#1}}

\def\received#1{\def\receiveddate{#1}}
\def\revised#1{\def\reviseddate{#1}}
\def\accepted#1{\def\accepteddate{#1}}

\def\asciiaddress#1{\def\theasciiaddress{#1}}

\long\def\asciiabstract#1{\long\def\theasciiabstract{#1}}

\let\\\par\let\thelognumber\relax\let\thevolumenumber\relax
\let\thepapernumber\relax\let\thevolumeyear\relax\let\startpage\relax
\let\finishpage\relax\let\publishdate\relax\let\receiveddate\relax
\let\reviseddate\relax\let\accepteddate\relax\let\theasciititle\relax
\let\theasciiauthors\relax\let\theasciiaddress\relax
\let\theasciiabstract\relax

\let\theasciiemail\relax

\ifplaintex
\font\logobig=cmssbx10 scaled 3836
\font\logomed=cmssbx10 scaled 2557
\else
\font\logobig=cmssbx10 scaled 4200
\font\logomed=cmssbx10 scaled 2800
\fi

\long\def\makeagttitle{   
\count0=\startpage
\agt\hfill      
\hbox to 45truept{\vbox to 0pt{\vglue -13truept{\logomed A\kern -.37em{\logobig 
T}\kern -.38em G}\vss}\hss}
\break
{\small Volume \thevolumenumber\ (\thevolumeyear)
\startpage--\finishpage\nl
Published: \publishdate}

\vglue .25truein

{\parskip=0pt\leftskip 0pt plus
1fil\def\\{\par\smallskip}{\Large\bf\thetitle}\par\medskip} \vglue
0.05truein

{\parskip=0pt\leftskip 0pt plus 1fil\def\\{\par}{\sc\theauthors}
\par\medskip}%
 
\vglue 0.03truein 

{\small\leftskip 25truept\rightskip 25truept{\bf Abstract}\stdspace\theabstract

{\bf AMS Classification}\stdspace\theprimaryclass
\ifx\thesecondaryclass\relax\else; \thesecondaryclass\fi\par
{\bf Keywords}\stdspace \thekeywords\par}\vglue 7truept

}   

\ifplaintex
\font\phead=cmsl9 scaled 950
\font\pnum=cmbx10 scaled 913
\font\pfoot=cmsl9 scaled 950
\headline{\vbox to 0pt{\vskip -4.5mm\line{\small\phead\ifnum
\count0=\startpage ISSN 1472-2739 (on-line) 1472-2747 (printed)
\hfill {\pnum\folio}\else\ifodd\count0\def\\{ }%
\ifx\theshorttitle\relax\thetitle\else\theshorttitle\fi\hfill{\pnum\folio}
\else\def\\{ and }{\pnum\folio}\hfill\ifx\theshortauthors\relax\theauthors
\else\theshortauthors\fi\fi\fi}\vss}}
\footline{\vbox to 0pt{\vglue 0mm\line{\small\pfoot\ifnum\count0=\startpage
\copyright\ \gtp\hfill\else
\agt, Volume \thevolumenumber\ (\thevolumeyear)\hfill\fi}\vss}}
\else
\font=cmsl9 scaled 1050
\font\lnum=cmbx10 
\font=cmsl9 scaled 1050
\makeatletter
\def\@oddhead{{\small\lhead\ifnum\count0=\startpage ISSN 1472-2739 
(on-line) 1472-2747 (printed)\hfill {\lnum\number\count0}\else\ifodd\count0
\def\\{ }\ifx\theshorttitle\relax \thetitle \else\theshorttitle\fi\hfill
{\lnum\number\count0}\else\def\\{ and }{\lnum\number\count0}
\hfill\ifx\theshortauthors\relax 
\theauthors\else\theshortauthors\fi\fi\fi}}\def\@evenhead{\@oddhead}
\def\@oddfoot{\small\lfoot\ifnum\count0=\startpage\copyright\ \gtp\hfill\else
\agt, Volume \thevolumenumber\ (\thevolumeyear)\hfill\fi}
\def\@evenfoot{\@oddfoot}
\makeatother
\fi
\let\maketitlepage\makeagttitle
\let\makeshorttitle\maketitlepage
\let\maketitle\maketitlepage

\newwrite\gtoutfile
\long\gdef\makeheadfile{  
{\def\\{, }\def\s{ }
\immediate\openout\gtoutfile head.xxx
\immediate\write\gtoutfile{Proxy-for: \ifx\theasciiauthors\relax
\theauthors\else\theasciiauthors\fi\s<\ifx\theasciiemail\relax\theemail\else\theasciiemail\fi>}
\immediate\write\gtoutfile{\noexpand\\}
\immediate\write\gtoutfile{Authors: \ifx\theasciiauthors\relax
\theauthors\else\theasciiauthors\fi}
{\def\\{ }\immediate\write\gtoutfile{Title: \ifx\theasciititle\relax
\thetitle\else\theasciititle\fi}}
\immediate\write\gtoutfile{Subj-class: GT or SG, GR etc}
\immediate\write\gtoutfile{MSC-class: \theprimaryclass\ifx\thesecondaryclass\relax\else, \thesecondaryclass\fi}
\immediate\write\gtoutfile{Journal-ref: Algebraic and Geometric Topology \thevolumenumber\s
(\thevolumeyear) \startpage-\finishpage}
\immediate\write\gtoutfile{Comments: Published by Algebraic and
Geometric Topology at}
\immediate\write\gtoutfile{\s\s\s  http://www.maths.warwick.ac.uk/agt/AGTVol\thevolumenumber/agt-\thevolumenumber-\thepapernumber.abs.html}
\immediate\write\gtoutfile{\noexpand\\}
\immediate\write\gtoutfile{}
\ifx\theasciiabstract\relax
\immediate\write\gtoutfile{\theabstract}\else
\immediate\write\gtoutfile{\theasciiabstract}\fi
\immediate\write\gtoutfile{}
\immediate\write\gtoutfile{\noexpand\\}
\immediate\write\gtoutfile{}
\immediate\closeout\gtoutfile}}  

\def\maketitlepage{\makeagttitle\makeheadfile}
\let\makeshorttitle\maketitlepage
\let\maketitle\maketitlepage

\lognumber{23}
\volumenumber{3}
\volumeyear{2003}
\papernumber{23}
\published{4 July 2003}
\pagenumbers{677}{707}
\received{7 October 2002}
\revised{2 May 2003}
\accepted{9 June 2003}

\usepackage{amssymb,graphicx}

\def\pv{\par\medskip}
\let\tm\rm
\def\d #1{\mathbb{#1}}

\def\r #1{\mathcal{#1}}

\def\ad #1{\hbox{$\overline #1$}}
\def\ch #1{\hbox{$\widehat #1$}}
\def\bo #1{\hbox{$\partial #1$}}

\font\ssmall cmr8

\def\cmp{com\-po\-nent}

\def\cs{Schar\-le\-mann cycle}

\def\res{res\-pec\-tive\-ly}

\def\nd{knot}

\def\Bequiv{$B_{/_{\sim_\partial}}$}
\let\cqfd\endproof

\newenvironment{Relax}{\relax}{\relax}
\begin{document}

\def\boxit#1#2{\vbox{\hrule\hbox{\vrule
\vbox spread#1{\vfil\hbox spread#1{\hfil#2\hfil}\vfil}%
\vrule}\hrule}}

\long\def\emboite#1#2{\setbox4=\vbox{\hsize #1cm\noindent \strut
#2}\boxit{4pt}{\box4}}

\let\dps\displaystyle

\def\stdfigure #1 by #2 (#3 scaled #4) (#5: #6)
{\begin{figure}[ht!]
\cl{\includegraphics[width=#1]{Figure#3.eps}}
\caption{#6}
\end{figure}}

\title{Thin presentation of knots and lens spaces}

\author{A. Deruelle\\D. Matignon}

\address{Universit\'e D'Aix-Marseille I,
C.M.I. 39, rue Joliot Curie\\
Marseille Cedex 13, France}
\asciiaddress{Universite D'Aix-Marseille I,
C.M.I. 39, rue Joliot Curie\\
Marseille Cedex 13, France}

\email{deruelle@cmi.univ-mrs.fr, matignon@cmi.univ-mrs.fr}

\begin{abstract}
This paper concerns thin presentations of knots $K$ in closed
$3$-manifolds $M^3$ which produce $S^3$ by Dehn surgery, for some
slope $\gamma$.  If $M$ does not have a lens space as a connected
summand, we first prove that all such thin presentations, with respect
to any spine of $M$ have only local maxima.  If $M$ is a lens space
and $K$ has an essential thin presentation with respect to a given
standard spine (of lens space $M$) with only local maxima, then we
show that $K$ is a $0$-bridge or $1$-bridge braid in $M$; furthermore,
we prove the minimal intersection between $K$ and such spines to be at
least three, and finally, if the core of the surgery $K_\gamma$ yields
$S^3$ by $r$-Dehn surgery, then we prove the following inequality:
$\vert r\vert\leq 2g$, where $g$ is the genus of $K_\gamma$.
\end{abstract}
\asciiabstract{
This paper concerns thin presentations of knots K in closed
3-manifolds M^3 which produce S^3 by Dehn surgery, for some
slope gamma.  If M does not have a lens space as a connected
summand, we first prove that all such thin presentations, with respect
to any spine of M have only local maxima.  If M is a lens space
and K has an essential thin presentation with respect to a given
standard spine (of lens space M) with only local maxima, then we
show that K is a 0-bridge or 1-bridge braid in M; furthermore,
we prove the minimal intersection between K and such spines to be at
least three, and finally, if the core of the surgery K_gamma yields
S^3 by r-Dehn surgery, then we prove the following inequality:
|r| <= 2g, where g is the genus of K_gamma.}

\keywords{Dehn surgery, lens space,
 thin presentation of knots, spines of 3-manifolds}

\primaryclass{57M25}
\secondaryclass{57N10, 57M15}
\maketitle

\begin{Relax}\end{Relax}

\section{Introduction}

All $3$-manifolds are assumed to be
compact, connected and orientable.
A {\it link} in a $3$-manifold is a compact and closed $1$-submanifold.
A {\it Dehn surgery} on a link $\cal L$ in a $3$-manifold $M$,
consists on removing a regular neighbourhood $N(\r L)$ of $\r L$, and
gluing back solid tori to the corresponding toroidal
boundary components of $M-N(\r L)$ by boundary-homeomorphisms.
In [27, 41] Wallace and Lickorish have proved independently that a
compact, connected and orientable $3$-manifold
can be obtained by Dehn surgery on a link in
the $3$-sphere $S^3$.
Dehn surgery on knots (one-component links)
are of high interest in low dimensional topology, see the nice
surveys of Gordon [19] or Luecke [28].

In this paper, we are interested in $3$-manifolds obtained
by Dehn surgery on a knot in $S^3$ and in particular,
in the following question:
what do the knots look like in an arbitrary closed $3$-manifold
if they produce $S^3$ by Dehn surgery?
We will answer the question towards the thin presentation of the knots
according to a spine of the $3$-manifold.

Each closed $3$-manifold $M$ is a $3$-ball with an identification on its
boundary (see Section 2 and [33, Chapter~2] for details).
Let $\Sigma$ be the corresponding spine of $M$;
i.e.\ the identified boundary of the $3$-ball.
Then  $M$ where $\Sigma$ and an interior point are removed,
is homeomorphic to $S^2\times \d R$.
We consider this $2$-spheres foliation of $M$,
and to study what the knots look like,
 we define their thin presentations in $M$,
in a similar way as Gabai did for knots in $S^3$~[15, Section~4.A],
but with respect to the spine $\Sigma$.

This notion is very useful, and has played a
key-point in the proof of the property $\r R$ by Gabai [15] , as well as
in the solution of the complement problem by Gordon and Luecke [22].
This concept has been used also in other
important $3$-dimensional topology problems,
as the recognition of $S^3$ by Thompson [39] or
the study of Heegaard diagram of the $I$-fibered on surfaces by
Scharlemann [37].  Now, the thin presentation of knots is in itself
the topic of many works (see for example [1, 24, 35, 40, 45]).

The first result of the present paper is the following.

\let\tmg\bf
{\tmg Theorem 1.1}\qua
{\sl Let $K$ be a knot in a closed $3$-manifold $M$,
that does not have a lens space as a connected summand.
If there exists a spine $\Sigma$
such that a thin presentation of $K$ in $M$, with respect to $\Sigma$,
has a local minimum then $K$ cannot yield $S^3$ by Dehn surgery.}

Let put this result in terms of knots in $S^3$,
giving this other formulation of Theorem~1.1.

\pv {\tmg Theorem 1.1.Bis}\qua
{\sl Let $k$ be a knot in $S^3$.
Let $k(\alpha)$ be the $3$-manifold obtained
by $\alpha$-Dehn surgery on $k$ and $k_\alpha$ be the core of the surgery.
If $k(\alpha)$ does not contain a lens space summand then,
for any spine of $k(\alpha)$, all the thin presentations of $k_\alpha$
have only local maxima.
}

\pv
Recall that for any knot in $S^3$,
only two slopes can produce a {\it reducible} manifold
(i.e.\ containing a $2$-sphere that does not bound a $3$-ball)
by [23] (see also [25] for an alternative proof);
and similarly for non-torus knots,
at most three slopes can produce a lens space [8, 32].

So, Theorem~1.1.Bis implies that
for the cores $K$ of all Dehn surgeries on a knot in $S^3$
but a finite number, their thin presentations, with respect to all spines
in the surgered manifold, have only local maxima.

Now, the main part of the paper is devoted to the case
where $M$ is a lens space $L$ within we define {\it standard} spines.
We know Dehn surgeries on the trivial knot to produce $S^3$,
$S^2\times S^1$ and general lens spaces.
So, the problem is focus on Dehn surgeries on non-trivial knots in $S^3$
and in particular, is it possible to obtain a lens space?
We know the answer to be negative for $L(1,m)=S^3$ [22],
and also for $L(0,1)=S^2\times S^1$ [15].
In the general case, the answer is positive for many knots [3, 16].
Nevertheless, the question whether Dehn surgery on a knot in $S^3$ produces
a lens space, is still open and
subject to a large sphere of investigations [14, 17, 19, 28].

The problem is completely solved for torus knots [32] and
satellite knots [6, 20, 42, 43].
It is also known that there are many hyperbolic knots
which produce lens spaces;
among them the $(-2,3,7)$-pretzel knot [14] produces $L(18,5)$ and
$L(19,7)$.
Furthermore, Berge in~[3] exhibits infinite families of knots with
a Dehn surgery yielding a lens space and gives its construction.
In~[19], Gordon asked Question~5.5:
Does every knot $K$ producing a lens space for some Dehn surgery
appear in Berge's list?
As there is no known example concerning the production of a lens space
with order smaller than five, an affirmative answer to this question
would imply the following conjecture to be true.

\pv
{\tmg Conjecture A}\qua (Gordon~'90 [19, Conjecture~5.6])

{\sl Dehn surgery on a non-trivial knot in $S^3$ cannot yield a lens space
with order less than five.}

\pv
A knot in a lens space $L$ is a $n$-bridge braid if,
for a Heegaard solid torus $V$ of $L$ (i.e.\ $L-V$ is a solid torus),
it can be isotoped to a braid in $V$ which lies in $\partial V$
except for $n$ bridges~[16].

Then a $0$-bridge braid is a torus knot (in $\partial V$).
And a knot is a $1$-bridge braid if it is the union of two arcs $\alpha$ and
$\beta$, each transverse to the meridional disks of $V$, such that:
$\alpha$ is lying on $\partial V$ and $\beta$ is properly embedded in $V$
and is cobounding a disk in $V$ with an arc on $\partial V$.

In [3], Berge asked a question about the production of lens spaces,
but in terms of a knot in the lens space:
If $k$ is a knot in a lens space such that Dehn surgery
on $k$ yields $S^3$, must $k$ be a $0$ or $1$-bridge knot in the lens space?

Let us remark that Berge also proves that a $1$-bridge knot in a lens space
(i.e.\ a $(1,1)$-knot),
producing $S^3$ by Dehn surgery is isotopic to a knot
which is simultanously braided
with respect to both of the solid tori of genus one Heegaard splitting
of the lens space.
Many works concern $(1,1)$-knots, see for example
[12, 13].

Following Berge and Gordon, one would state the following conjecture
which places the point of view in terms of knots in lens spaces.
\pv
{\tmg Conjecture B}
\qua
{\sl If a knot $K$ in a lens space $L$ produces $S^3$ by Dehn surgery then
$K$ is a $0$ or $1$-bridge braid.
}
\pv
In this framework, where $M=L$ is a lens space,
we define a {\it thin presentation} with respect to a standard spine
$\Sigma$ of $L$.
Then a local maximum in a thin presentation
is {\it inessential} if one can isotope it to $\Sigma$.
And it is {\it essential} if it cannot be isotoped.
After what, we introduce an {\it essential} thin presentation based
on the existence of such essential local maxima in the first ones.
For more details, we refer to Sections~2 and 4.

Note that in all the following, we consider $L$ different from $S^3$ and
also from $S^2\times S^1$.
We prove the following result.
\pv
{\tmg Theorem 1.2}\qua
{\sl
Let $K$ be a knot in $L$ yielding $S^3$ by Dehn surgery.
If there exist a standard spine $\Sigma$ and
an essential thin presentation of $K$ with respect to $\Sigma$
beginning by a local maximum,
then $K$ is a $0$ or $1$-bridge braid in $L$.}
\pv
Let $K$ be a knot in a lens space $L$ and
$\Sigma$ be a standard spine of $L$.
If a thin presentation of $K$ with  respect to $\Sigma$
has only inessential local maxima then
$K$ can be isotoped onto $\Sigma$;
the authors refer again to Section~4 for the definition.
For convenience, we say that $K$ is a {\it standardly spinal knot}
in $L$.

So, in the light of Theorem~1.1, we state the following conjecture.
\pv
{\tmg Conjecture C}
\qua
{\sl If $K$ is a knot in a lens space $L$ yielding $S^3$ by Dehn surgery
then $K$ is a standardly spinal knot.
}
\pv
{\tmg Q{uestion} D}
\qua
{\sl If $K$ is a standardly spinal knot,
must $K$ be a $(1,1)$-knot?
}
\pv
The {\it $\d RP^3$--Conjecture}  (i.e.\ Conjecture~A
for real projective $3$-space) claims that
$\d RP^3$ cannot be obtained by Dehn surgery on a non-trivial knot in $S^3$.
Let us note here that if one can prove Conjecture~C and
answers positively to Question~D, then
it would imply Conjecture~B and so the $\d RP^3$-- Conjecture.

Let $s$ be the minimal geometric intersection number between $\Sigma$ and
$K$.
Note that in [9], where $L=L(2,1)=\d R P^3$, it is shown that
if a thin presentation of $K$ with respect to a minimal projective plane
(as standard spine) has only local maxima then $s=1$ and therefore,
the core of the surgery is the trivial knot in $S^3$.
This result can now be viewed as a consequence of Theorem~1.2.

We know the $\d RP^3$-conjecture to be satisfied for cable knots [42, 43].
Furthermore, the standard spine of $\d RP^3$ is a projective plane and
by [11], we know $s\geq 5$.

In this paper, we also look at the number of intersections $s$,
but for a knot in a general lens space.
\pv
{\tmg P{roposition} 1.3}
\qua
{\sl
If $K$ is neither a $0$ nor $1$-bridge braid in $L$,
then $s\geq 3$ for all standard spines of $L$.
}

\pv
If the core of the surgery $K_\gamma$ is not a torus knot in $S^3$,
the slope $r$ that yields the lens space $L$ is an integer [8].
Let $g$ be the genus of $K_\gamma$.
In [17], Goda and Teragaito show that $\vert r\vert\leq 12g-6$,
if $K_\gamma$ is hyperbolic, and conjectures that $\mid r\mid\leq 4g-1$.
This inequality has recently been improved by Ichihara~[26]:
$\vert r \vert\leq 3\cdot 2^{7\over4} g$.
Here, we prove an inequality involving also the genus $g$ and
the slope $r$ but towards non $0$ nor $1$-bridge braids.

Let us mention that if $K_\gamma$ is a torus knot then
$K_\gamma$ is a $0$-bridge braid in $S^3$ and
so $K$ is a $0$-bridge knot in $L$ by [16].

\pv
{\tmg T{heorem} 1.4}
\qua
{\sl If $K$ is not a standardly spinal knot in $L$,
then $\mid r\mid\leq 2g$.
}

\pv
The main results of this paper are based on intersection graphs techniques
[8, 22] and Cerf Theory, in a similar way as Gordon and Luecke~[22],
proving that knots in $S^3$ are determined by their complements.
Let give a brief description of these arguments.

Let $K$ be a knot in a closed $3$-manifold $M$,
which produces $S^3$ by Dehn surgery.
We define a {\it spinal presentation} of the $3$-manifold $M$
which allows us to define a {\it thin presentation} of knots in $M$.
Therefore, we obtain, on one side, a $M$-foliation
(with level $2$-spheres, according to a height function)
in which $K$ is in thin presentation,
and on the other side, a $S^3$-foliation
in which $K_\gamma$ is in a thin presentation, by~[15].

Then,
we study the intersection of two one-parameter families of surfaces
whose we deduce the respective foliations,
to find a pair of properly embedded surfaces in
the complement of $N(K)$ in $M$
(where $N(K)$ is a regular neighbourhood of $K$).
This pair of surfaces gives rise to a pair of intersection graphs,
in the usual way [7, 16, 22].
A study of the two foliations leads to properties for the associated graphs.
And conversely,
a study of the pair of intersection graphs leads to
some properties of the corresponding foliations.
Comparing these properties with the original gives a contradiction or
the required properties of the knot ($0$ or $1$-bridge braid).

The paper is organized as follows.
In Section~2,
we define the basic tools of the paper:
the {\it thin presentation} of knots in closed $3$-manifolds
associated to a particular foliation and
the corresponding {\it essential} thin presentation in the case of
lens spaces;
also the {\it intersection graphs} and
the links between foliations and
properties of the graphs.
Section~3 is devoted to the proof of Theorem 1.1.

In the last sections,
we only consider the case where $M$ is a lens space $L$ say.
In Section~4,
we prove that if a thin presentation of $K$ in $L$,
with respect to a standard spine $\Sigma$,
begins by an {\it essential} local maximum then
$K$ intersects $\Sigma$ only once.
In Section~5,
we first extend the result of the previous section proving as a consequence,
Theorem~1.2. Then, as
a converse and then,
focusing on the number of intersections
between $K$ and the standard spines in $L$, we prove Proposition~1.3.
Finally, in Section~6, we use the previous results to prove Theorem~1.4.

\pv
{\tmg Acknowledgement}\qua
We address our deep thanks to Mario Eudave-Mu\~noz for interesting and
helpful discussions.

\section{Preliminaries}

In this section, we define the
common background and fix the notations for all the following sections.

Let us first recall the definition of Dehn surgery.

If $k$ is a  knot in $S^3$,
we denote $X_k=S^3-$int$N(k)$ the exterior of the knot
(also called the {\it space of the knot}),
where $N(k)$ is a regular neighbourhood of $k$.
So, the boundary of $X_k$ is a torus $T_k$ and
a {\it slope} $r$ on $T_k$ is the isotopy class of an un-oriented essential
simple closed curve.
The slopes are then
parametrized by $\d Q\cup \{\infty\}$ (for more details, see [36]).

A {\it $r$-Dehn surgery} on $k$ consists in gluing a solid torus
$V=S^1\times D^2$ to $X_k$ along $T_k$ such that
$r$ bounds a meridional disk in $V$.
We denote $X_k(r)$ the resulting closed $3$-manifold.
The core of $V$ becomes a knot $k_r$ in $X_k(r)$
called the {\it core of the surgery}.

Dehn surgery on a knot in a closed $3$-manifold is defined in a similar way,
 by gluing a solid torus to the exterior of the knot such that
the chosen slope bounds a meridional disk (in the attached solid torus).
Note that if we do \hbox{$r$-Dehn} surgery on a knot $k$ in $S^3$,
then we can obtain $S^3$
by doing Dehn surgery on $k_r$ in the closed $3$-manifold $X_k(r)$.

\sh{Thin presentation of knots in $S^3$}

For convenience, we recall
the definition of a thin presentation of knots in $S^3$,
introduced by Gabai [15, Section~4.A].

\begin{figure}[ht!]
\cl{\rotatebox{90}{\includegraphics[width=50mm]{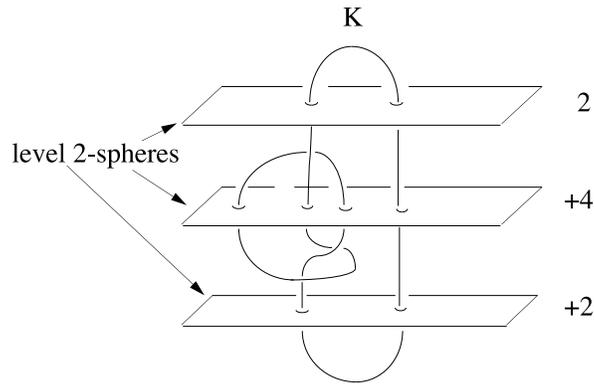}}}
\caption{Thin presentation of complexity $8$}
\end{figure}

 If $\pm\infty$ are the North and South poles of $S^3$, note that
$S^3-\{\pm\infty\}\cong S^2\times\d R$.
Then we have an associated height function $h:
S^3-\{\pm\infty\}\longrightarrow \d R$
which is the projection onto the second factor.
A sphere $\ch P_t=S^2\times\{ t\}$ in such a foliation for $t\in \d R$,
is called a {\it level $2$-sphere}.

Let $k$ be a knot in $S^3$.
By an isotopy of $k$, we may assume that
$k\subset S^3-\{\pm\infty\}$ and that
$h\vert_k$ is a Morse function, that is,
$h\vert_k$ has only finitely many critical points,
all non-degenerate, with all critical values distinct.
Each critical value represents a tangency point between
the corresponding level $2$-sphere and the knot.

Between each pair of consecutive critical values of
$h{\vert}_k$, the level $2$-spheres have the same
geometric intersection number with the knot.
Given such a Morse presentation of $k$, let $S_1,\dots,S_m$ be level
$2$-spheres, one between each pair of consecutive critical levels.

One then calls the number
$\displaystyle\sum_{i=1}^{m}\vert S_i\cap k\vert$
the {\it complexity} of the Morse presentation.
A {\it thin presentation} of $k$ is
a Morse presentation of minimal complexity (Figure 1).

A properly embedded surface in $X_k$,
isotopic to $P_t=\ch P_t\cap X_k$ is called
a {\it level surface} of the presentation, and
$\partial P_t$ consists of several parallel copies
of a meridional curve on $\partial X_k=\partial N(k)$.

\sh{Spinal presentation of closed $3$-manifolds}

A classical method for constructing any closed $3$-manifold consists
in matching up all the $2$-simplices in the boundary
of a triangulated $3$-cell;
known as the {\it maximal cave} method~[33, Chapter~2].

Let $M$ be a closed $3$-manifold.
Then we can see $M$ as \Bequiv,
where $B$ is a (closed) $3$-ball and
$\sim_\partial$ is an equivalence relation defined on
the $2$-sphere $S=\partial B$.
We call $\Sigma=S_{/\sim_\partial}$ the {\it spine} of $M$.
Let recall the construction of $M$ as \Bequiv;
for more details, see [29, 33].

All $3$-manifolds are triangulable~[5, 31].
So, let $\r T$ be a triangulation of $M$.
If $N$ is a combinatorial sub-manifold of $M$,
we denote $\mid N\mid$ the sub-complex of $\r T$,
corresponding to the closure of $N$ in $M$.
Let $B_0$ be the interior of a $3$-simplex $\sigma_0$ in $\r T$.
Let us choose an open $3$-simplex $\sigma_1$ in $\r T-B_0$ such that
$\mid B_0\mid\cap\ \ad\sigma_{1}$ is a (non-empty) union of
$2$-simplices in $\partial\sigma_1$.
And set $B_1$ to be the union of $B_0$ and $\sigma_1$, glued along
the interior of one of the simplices of $\mid B_0\mid\cap\ \ad\sigma_{1}$
(just choose one).
Now extend this construction by induction in the following way:
For each integer $k$,
let $B_{k+1}$ be the union of $B_k$ and
an open $3$-simplex $\sigma_{k+1}$ in $\r T-B_k$
as described above;
that is, one have to choose an open $3$-simplex
$\sigma_{k+1}$ in $\r T-B_k$ and his ``prefered'' $2$-simplex in
$\mid B_k\mid\cap\ \ad\sigma_{k+1}$.
Then, for each integer $k$, $B_k$ is opened and
$\vert B_k\vert$ is a closed $3$-ball in $M$.

By this process,
we must include all the $3$-simplices of $\r T$,
because $M$ is connected.
Furthermore, $M$ is compact so,
there exists an integer $N$ such that
$M-B_N$ contains no $3$-simplex of $\r T$.
If $B$ is the closure of $B_N$,
then $B$ is a closed triangulated $3$-ball
(with the induced triangulation on the boundary).
But note that if $M-B_N$ is the union of the $2$-simplices of $\partial B$
represented only once,
these $2$-simplices are represented twice on $\partial B$.
So, they define an equivalence relation on $\partial B$,
we denote $\sim_\partial$.
Furthermore,
we then have $M=\hbox{\Bequiv}$, and $\partial B$ with
the identified $2$-simplices is exactly the spine $\Sigma=M-B_N$.

We say that $\Sigma$ is {\it canonical},
by meaning that $S$ is triangulated, and
the $2$-simplices are identified by the equivalence relation.
All spines are assumed to be canonical in the following.

 In the case of lens spaces,
we say that a spine is {\it standard} if
it corresponds to the usual, refering to Rolfsen~[36, p.236].
That is,
each $2$-simplex in the $2$-sphere $S$
has one edge on the equator circle of $S$,
and one vertex at either the North or South pole.
Moreover,
each simplex in the north hemisphere of $S$
is identified to a simplex in the south one,
by some ${2\pi q\over p}$-rotation~(Figure~2).
The edges on the equator are all identified
to a single embedded circle $\kappa$ in $\Sigma$.
We call $\kappa$, the {\it core} of the standard spine and
this is the singular set of $\Sigma$.
We would like to warn the reader on
the special notation used for the core of a standard spine:
we use the greek letter kappa $\kappa$;
not to misunderstand with the notation of knot $K$.

A regular neighbourhood of $\kappa$ in $\Sigma$ is defined to be a
{\it spinal helix}:
this is a $2$-complex obtained by removing a disk in $\Sigma$
disjoint from $\kappa$.
In $L(2,1)$, it is a M\"obius band.
But in the general case $L(p,q)$ for $p>2$,
 a spinal helix is a $2$-complex with a singular set $\kappa$ on which
the ``surface'' runs a
finite number of times;
that what we call the {\it order} of the spinal helix.
If it is obtained from a spine, by removing a disk,
the spinal helix is of order $p$,
which is the order of the spine and also of the lens space $L(p,q)$.

\stdfigure 35mm by 273mm (02 scaled 200) (Figure~2: Standard spine of
{$L(p,q)$} and core $\kappa$)

\sh{Thin presentation of knots in closed $3$-manifolds}

Now, we define the thin presentation of a knot in $M$
with respect to a spine $\Sigma$.
We note $M=\hbox{\Bequiv}$ according to the spine
$\Sigma=\partial B_{/\sim_\partial}$.

Let $\infty$ be an interior point of the $3$-ball $B$.
Then $M-(\Sigma \cup\{\infty\})\cong S^2\times\d R_{>0}$.
So, we have an associated
height function $h: M-(\Sigma\cup \{\infty\})\rightarrow \d R_{>0}$
which is the projection on the second factor.
We extend $h$ to $\Sigma$ by setting $h\vert_{\Sigma}=0\in\d R_{\geq0}$,
so $h(M-\{\infty\})=\d R_{\geq0}$.
The {\it level $2$-spheres} are the spheres
$\ch Q_t=S^2\times\{ t\},\ t\in ]0,+\infty[$.

Let $K$ be a knot in $M$. By an isotopy of $K$, we may assume that
$K\subset M-\{\infty\}$ is transverse to $\Sigma$ and $h\vert_{K}$
is a Morse function.

Similarly as for knots in $S^3$, between each pair of consecutive critical
values of $h\vert_{K}$,
the level $2$-spheres of the foliation have the same geometric intersection
number with the knot (Figure~3).
Given such a Morse presentation of $K$,
let $S_1,\dots,S_m$ be level $2$-spheres,
one between each pair of consecutive critical levels.
Furthermore, let $S_0$ be a level $2$-sphere between $\Sigma$ and
the first critical level.
So, $\widehat Q_0=S_0=\partial N(\Sigma)$, where $N(\Sigma)$ is
a regular neighbourhood of $\Sigma$.

\stdfigure 100mm by 96.6mm (03 scaled 398) (Figure~3:
Morse presentation in a $3$-manifold)

We then call the number
$\displaystyle\#(K)=\sum_{i=0}^{m}\vert S_i\cap K\vert$
the {\it complexity} of the Morse presentation.
A {\it thin presentation} of $K$ is a
Morse presentation of minimal complexity;
that is, precisely a presentation of the knot obtained
by further isotopies of $K$
in $M-\Sigma\cong M-IntN(\Sigma)$
on a given Morse presentation to minimize the complexity.
In a thin presentation, one cannot decrease
$\displaystyle\#(K)-s_0=\sum_{i=1}^{m}\vert S_i\cap K\vert$.

We denote $M_K=M-$int$N(K)$, the exterior of the knot $K$.
A properly embedded surface in $M_K$, isotopic to $Q_t=\ch Q_t\cap M_K$ is
called a {\it level surface} of the presentation.
Remark that $Q_0=\ch Q_0\cap M_K=S_0\cap M_K$ is
(still) in a regular neighbourhood of $\Sigma\cap M_K$.

Note first that $\Sigma$ is not necessary a minimal spine;
i.e.\ a spine with a minimal intersection number with the knot
among all spines of the $3$-manifold.
Thin presentations of $K$ are defined with an arbitrary choosen spine
$\Sigma$.

And finally note that if, for a presentation,
$K$ intersects the singular set of $\Sigma$ then
$h\vert_{ K}$ is not a Morse function
(because of the openess property of being a Morse function).

We will see, in Section~4 that in certain conditions
we can minimize the intersection between $K$ and
a standard spine $\Sigma$ of a lens space,
by allowing intersection with $\kappa$, the singular set of $\Sigma$.

\sh{Associated intersection graphs}

Let $k$ be a knot in $S^3$ and $M=X_k(r)$ the $3$-manifold
obtained by $r$-Dehn surgery on $k$. Recall that $k_r$ denotes the
core of the surgery and $X_k=S^3-$int$N(k)\cong M-$int$N(k_r)=M_{k_r}$.

Let us consider a thin presentation of $k$ in $S^3$ and
a thin presentation of $k_r$ in $M$, associated to any spine.

Let \ch P and \ch Q denote level $2$-spheres in the foliations of $S^3$ and
$M$, respectively and
$P=\ch P\cap X_k$ and $Q=\ch Q\cap X_k$ be the corresponding level surfaces.
Then $P$ and $Q$ are planar surfaces, properly embedded in $X_k$.
The torus boundary $T_k=\partial X_k$ contains the slopes $r$ and $\infty$,
and each component of $\partial P$ (resp.\ $\partial Q$)
represents $\infty$ (resp.\ $r$).
Moreover, up to isotopies,
$P$ and $Q$ are transverse, and
each component of $\partial P$ intersects each component of $\partial Q$
exactly $\Delta$ times,
where $\Delta$ is $\Delta(\infty,r)$,
the geometric intersection number between $r$ and $\infty$ on $T_k$.

Let us recall the construction of intersection graphs
coming from a pair of planar surfaces properly embedded in $X_k$;
this is described in details in [18, 22].

 Let $(G,H)$ denote the pair of graphs associated to $(\ch P,\ch Q)$.
Capping off the boundary components of $P$ (resp.\ $Q$) with
meridional disks of $N(k)$ (resp.\ of $N(k_r)$),
we regard these disks as defining the ``fat'' {\it vertices} of the graph
$G$ (resp.\ $H$) in \ch P (resp.\ in \ch Q).
The {\it edges} of $G$ (resp.\ $H$)
are the arc-components of $P\cap Q$ in $P$ (resp.\ in $Q$).

The endpoints of edges in $G$ and $H$ can be labelled in the following way.
We first number the \cmp s of $\bo P $ and $\bo Q$ in the order they appear
(successively) on $T_k$.
Let number the components of $\partial P$:
$V_1,V_2,\dots,V_p$;
and those of $\partial Q$:
$W_1,W_2,\dots,W_q$.
We then label the endpoints of an arc of $P\cap Q$ in $P$
(resp.\ in $Q$)
with the numbers of the corresponding \cmp s of $\partial Q$
(resp.\ of $\partial P$)
that intersect $P$
(resp.\ in $Q$)
to create these endpoints.
Thus, around each component of $\partial P$
(resp.\ $\partial Q$),
we see the labels $\{1,2,\dots,q\}$
(resp.\ $\{1,2,\dots,p\}$)
appearing in this order (either clockwise or anticlockwise).
So these labels of the arcs of $P\cap Q$
allows us to label the endpoints of edges in $G$ and $H$
whether the arcs are viewed in $P$ or $Q$, respectively.

A vertex is {\it positive} if the labels appear clockwise around it;
otherwise, we say it is {\it negative}.
And two vertices are {\it parallel} if they have the same {\it sign},
i.e.\ they are both positive or both negative;
otherwise they are {\it antiparallel}.

In this framework,
we say that the pair of graphs $(G,H)$ is of {\it type} $(\ch P,\ch Q)$.

If $P$ and $Q$ are orientable surfaces,
we have the so called {\it parity rule}~[8]:
an edge $e$ in $G$ joins two parallel vertices if and only if
$e$ joins two antiparallel vertices in H.

Let $G^0$ and $G^1$ denote the set of vertices and edges of $G$, \res.
A {\it face} of $G$ is the closure of
a connected \cmp\ of $\ch P -(G^1\cup G^0)$.
Similarly, we denote $H^0$ and $H^1$,
the set of  vertices and edges of $H$, respectively,
and we have the same definitions for the faces of $H$.

Let $\Gamma$ be $G$ or $H$.
{\it A cycle} in $\Gamma$ is a subgraph of $\Gamma$
homeomorphic to a circle when vertices are considered as points;
the {\it length} of the cycle is the number of its edges.
Note that if a cycle in $\Gamma$ bounds a face of $\Gamma$,
then this face is necessarily a disk;
we will then say that it is a {\it disk-face}.
Two edges are said to be {\it parallel} in $\Gamma$ if
they form a cycle of length two which bounds a disk-face of $\Gamma$.

Two particular cycles play a key-role in the following.
A {\it trivial loop} is a cycle of length one which bounds a disk-face;
see Figure~4(b).
A {\it \cs} is a cycle which bounds a disk-face,
such that for an orientation of the cycle,
all the edges have the same label at their sink, $x$ say,
and also at their source, $y$ say.
Consequently $y=x\pm 1 \ (mod\ \zeta)$
(where $\zeta$ is $p$ or $q$ according to $\Gamma$ is $H$ or
$G$ respectively).
These labels $\{x,y\}$ are called the {\it labels} of the \cs;
see Figure~5(b).
And we then say that this is a $\{x,y\}$-\cs.

\sh{Links between foliations and intersection graphs}

For a $M$-foliation, a level $2$-sphere \ch Q separates $M$
in two connected components;
one of these contains the spine $\Sigma$ and
we say it {\it below} \ch Q (or {\it below} the level of \ch Q),
the other component setting to be {\it above} \ch Q.
With these definitions, we have also implicitly set what do we mean
by a presentation {\it beginning} by a local maximum (resp.\ minimum).

We set {\it above} a level $2$-sphere \ch P in a $S^3$-foliation,
the component of $S^3-$\ch P containing $+\infty$;
if containing $-\infty$, it is {\it below} \ch P.

In this paragraph, we keep the hypothesis and notations of the previous.
Let $\{X, Y\}=\{P,Q\}$.

\begin{figure}[ht!]
\cl{\includegraphics[width=90mm]{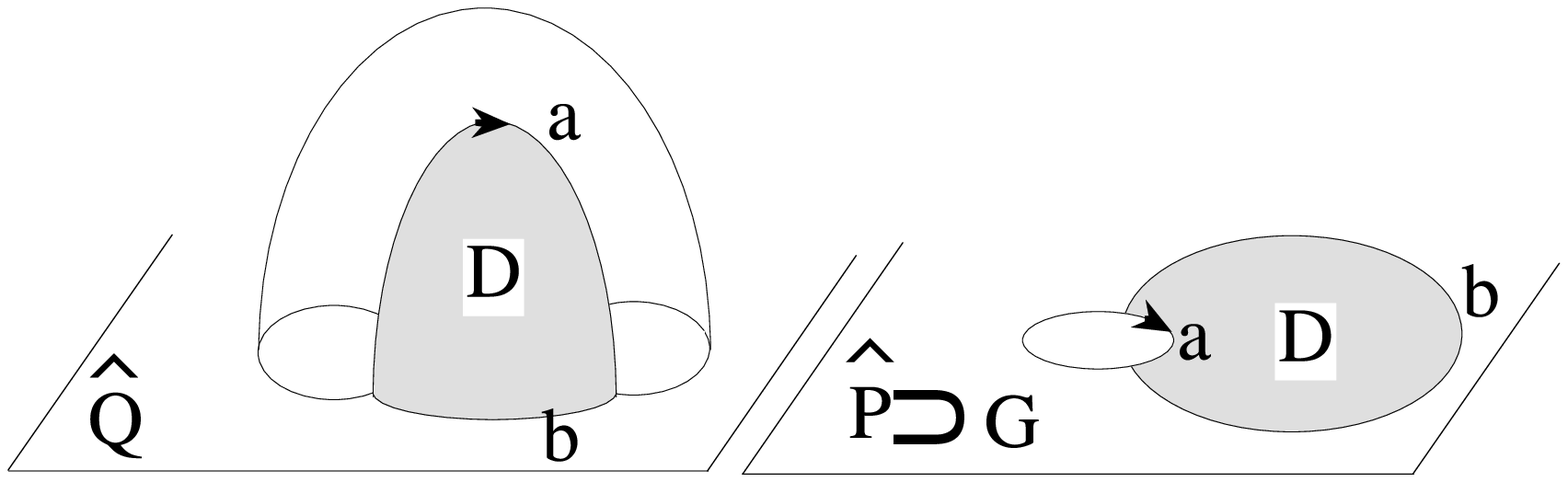}}
\medskip
\small\qquad\qquad\qquad (a)\qua High disk in $P$ \qquad\qquad\qquad
(b)\qua Trivial loop in $G$
\nocolon\caption{}
\end{figure}

Then, $X$ is said to be {\it high} (resp.\ {\it low})
with respect to $Y$ if
there exists a disk $D\subset Y$ below (resp.\ above) $X$ with
$\partial D=a\cup b$ such that $a=D\cap\partial X_k$ and
$b=D\cap X$ are simple arcs (Figure 4(a)).

The existence of such a (high or low) disk in $P$ (resp.\ in $Q$)
is equivalent to the existence of a trivial loop in the graph $G$
(resp.\ $H$);
see Figure~4(b).

Now, we define a characteristic corresponding to the existence of a
\cs.
We say that a disk $D\subset Y$ is {\it carrying} if
the two following conditions are satisfied,
for an orientation of $\partial D$ (Figure 5(a)).

\begin{figure}[ht!]
\cl{\includegraphics[width=100mm]{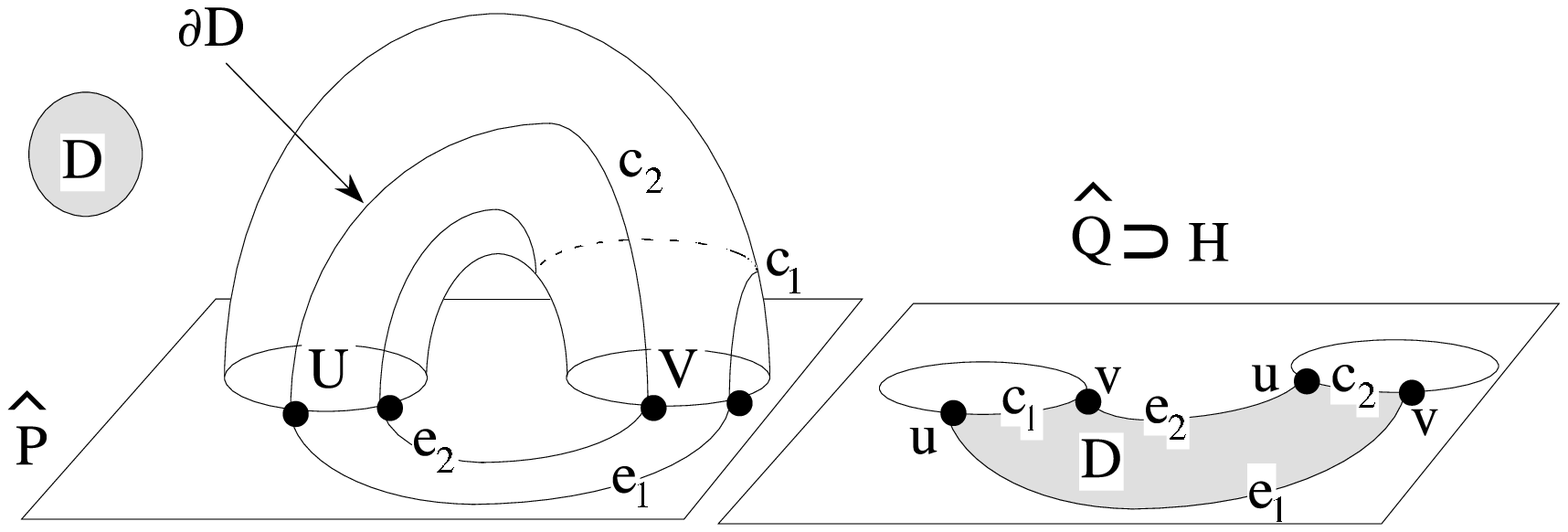}}
\medskip
\small\qquad\qquad\quad (a)\qua Carrying disk in $Q$\qquad\qquad
(b)\qua \cs\ in $H$
\nocolon\caption{}
\end{figure}

\items
\item[(i)] $D\cap X=\partial D\cap X$ is the disjoint union of
$n$ simple arcs, $e_1, \dots, e_n$
which join the same components $U$ and $V$ of $\partial X$ and
are all oriented from $V$ to $U$;

\item[\rm(ii)] $\partial D\cap \partial X_k$ is the disjoint union of
$n$ simple arcs, $c_1, \dots, c_n$,
all oriented from $U$ to $V$ in $\partial X_k$;
\enditems

We say that $X$ is  {\it carrier} with respect to $Y$ if
there exists a carrying disk $D\subset Y$.

The existence of such a carrying disk in $Q$ (resp.\ in $P$) is
equivalent to the
existence of a \cs\ in the graph $H$ (resp.\ $G$);
see Figure~5(b).
Actually, the vertices $U$ and $V$ in the above definition are
consecutive ones (in terms of labels) in the graph lying on $X$.

\section{The generic case}

Let $M$ be a closed $3$-manifold that does not contain a lens space as a
connected summand;
here the connected sum can be trivial,
that is we exclude the cases
where $M$ would be a lens space (or $S^3$ or $S^1\times S^2$).
Let $K$ be a knot in $M$ yielding $S^3$ by $\gamma$-Dehn surgery;
note $K_\gamma$ the core of the surgery.

\rk{Proof of Theorem 1.1}
Suppose for a contradiction that a thin presentation of $K$ in $M$
has a local minimum.

A {\it middle slab} in a thin presentation
of a knot is a family of level $2$-spheres, between consecutive
local minimum and local maximum
(of the thin presentation).
Let $K_\gamma$ be in a thin presentation in $S^3$.

Consequently,
the above thin presentations gave both  middle slabs,
denoted $\r M$ and $\r M'$ respectively.
Therefore, by [22, Proposition~1], there is a pair of level
$2$-spheres  $(\ch P,\ch Q)$ in $(\r M, \r M')$ such that $P$
is neither low
nor high with respect to $Q$, and vice-versa.
Equivalently, the pair of graphs $(G,H)$ of type $(\ch P,\ch Q)$,
 does not contain a trivial loop.
Note that \ch P and \ch Q are separating $2$-spheres because
they are level spheres.
Thus, we can apply the following combinatorial result, due to
Gordon and Luecke.

\pv
{\tmg T{heorem} 3.1}\qua [22, Proposition 2.0.1]
\qua
{\sl Let $(G,H)$ be a pair of intersection graphs of type $(S^2,S^2)$
without trivial loop.
If $H$ does not represent all types then $G$ contains a \cs, and
vice-versa.}
\pv
Moreover, if $H$ represents all types then $S^3$ contains a
$3$-sub-manifold
with non-trivial torsion ([22, Section 3] or [34]) which is impossible.
Consequently, $G$ contains a \cs, and so $M$ contains a lens space
as a connected summand by
[8, Lemma~2.5.2.b], which is the required contradiction proving
Theorem~1.1.

\sh{Intersection with all spines}

To conclude this section, we prove the following.

\pv
{\tmg L{emma} 3.2}
\qua
{\sl $K$ intersects all the spines of $M$.}
\pv
\proof
Let $M_K=M-$Int$N(K)$ and $M_K(\gamma)$ be the
$3$-manifold obtained by $\gamma$-Dehn surgery on $K$.
Assume $M_K(\gamma)=S^3$ and
note that $\gamma$ is not the meridional slope of $K$,
otherwise $M_K(\gamma)=M=S^3$.

We consider $M$ as \Bequiv, and
denote $\Sigma=S_{/\sim_\partial}$, the corresponding spine.
Let $s$ be the geometric intersection number between $K$ and $\Sigma$.
If $s=0$ then $K$ lies in Int$B$, and $M_K(\gamma)=M\#B_K(\gamma)$.
Consequently, $M_K(\gamma)$ is $M$ or
a reducible $3$-manifold, in contradiction with $S^3$.
\cqfd

In all the following, we consider the case where $M$ is a lens space.

\section{Lens space case}

Let $K$ be a knot in a lens space $L$,
and $\Sigma$ be a standard spine in $L$.
Let us suppose that $K$ is in a thin presentation,
with respect to $\Sigma$.
Let $\kappa$ be the core of $\Sigma$, that is the $1$-dimensional singular
sub-complex of $\Sigma$ (Figure~2).

Let $\mu$ be the level of a local maximum.
Denote $\alpha_\mu$
the arc on $K$ realizing this local maximum;
that is the arc on $K$ that
starts on the spine $\Sigma$, goes up straight through the level spheres,
passes tangently by the level $\mu$ and
goes down straight through the level spheres, back to $\Sigma$ (Figure~6).
Let $D_\mu$ be a disk properly embedded in $L_K-\Sigma$ where
$\partial D_\mu=\alpha\cup\beta$ is such that:

$\alpha=\partial L_K\cap D_\mu=\partial L_K\cap\partial D_\mu$ is an arc
parallel to $\alpha_\mu$,

$\beta=\Sigma\cap D_\mu=\Sigma\cap\partial D_\mu$ is an arc.

\noindent Remark that $\beta$ must intersect
the core $\kappa$ of $\Sigma$ in a several finite number of points.

\stdfigure 55mm by 145mm (06 scaled 280) (Figure~6 : Inessential local
maximum)

If there is such a disk then the corresponding local maximum at level $\mu$
is said to be {\it inessential}.
This means that the local maximum $\alpha_\mu$ can be isotoped,
to $\beta$ in
$\Sigma$.
For such an inessential local maximum $\mu_1$, let us do the isotopy of
$\alpha_{\mu_1}$ to $\Sigma$.
Then if there is another inessential local maximum $\mu_2$,
let us do the same.
And so on, until there is no more disk $D_\mu$ as we described above.
We then obtain a particular presentation of the knot $K$ that we call
{\it essential thin presentation}.

If an essential thin presentation of $K$ begins by a local maximum
then this local maximum is not inessential in the thin presentation;
i.e.\ it cannot be isotoped onto $\Sigma$.
We call it {\it essential} local maximum.

Let us note that an inessential local maximum is necessarily below
the first local minimum (if there is) in the thin presentation of $K$.
And also note that any local maximum above a local minimum must be
essential because of the thinness of the presentation of $K$.

Now, suppose that
$K$ is a knot in a lens space $L$ yielding $S^3$ by Dehn surgery.
The remaining of this section is devoted to the proof of the following
result.

\pv {\tmg T{heorem} 4.1}
\qua
{\sl If an essential thin presentation of $K$, with respect to the standard
 spine
$\Sigma$, begins by a local maximum
then $s=1$.}

\pv As usual, let $K_\gamma$ denote the core of the surgery in $S^3$ and
remark that $L_K=L-$Int$N(K)\cong S^3-intN(K_\gamma)=X_{K_\gamma}$.
So, let us consider a $S^3$-foliation in which $K_\gamma$ is in
thin presentation
and the $L$-foliation associated to an essential thin presentation of $K$
beginning by a local maximum, by hypothesis.
Let $\mu'$ be the first local maximum level of
this essential thin presentation of $K$.

Our goal is now to find two level surfaces $P=\ch P\cap L_K$, and
$Q=\ch Q\cap L_K$
in the $S^3$-foliation and the $L$-foliation
respectively, neither high, nor low, nor carrier, one with respect to
the other, and vice-versa.
Such a result is then in contradiction with Theorem 3.1.

To find this pair of transverse planar surfaces, we use the
Cerf Theory [7, Chapter 2] in a similar way as [22].

\sh{One-parameter families of $2$-spheres}

Recall that a middle slab in a thin
presentation of a knot is a family of level $2$-spheres,
between consecutive local minimum and local maximum.
We then consider:

\items
\item[\rm(i)] a middle slab $\{\ch P_\lambda\}_{\lambda\in I}$
in the $S^3$-foliation and

\item[\rm(ii)] a family of level $2$-spheres $\{\ch Q_\mu\}_{\mu\in J}$
in the
$L$-foliation between $\Sigma$ and the first local maximum.
\enditems

 We may suppose that $I=J=[0,1]$ and for convenience, we fix the index
notations $\lambda\in[0,1]$ for $S^3$ and $\mu\in[0,1]$ for $L$.

We denote a level $2$-sphere of a foliation $\cal H$, $\cal L$ or
$\cal C$ according to
its characteristic is $\cal H$igh, $\cal L$ow or $\cal C$arrier and
$\cal N$
if it is $\cal N$one of these.

\pv
{\tmg L{emma} 4.2}\qua
\items\sl
\item[\rm(i)] $Q_\mu$ cannot be $\cal L$ for all $\mu<\mu'$.

\item[\rm(ii)] There exists $\mu_0<\mu'$ such that $Q_\mu$ is $\cal C$ for all
$\mu<\mu_0$.
\enditems

\proof
(i)\qua
This follows immediatly from the fact
that the local maxima are essential in the thin presentation.

(ii)\qua
By the previous point, there exists $\mu_0<\mu'$ such that $Q_\mu$
is not $\cal L$ with respect to
$P_\lambda$, for all $\mu<\mu_0$ and all $\lambda\in[0,1]$.
Then  $Q_\mu$ is neither $\r L$ nor $\r H$, for all $\mu<\mu_0$.
Therefore, by [15, Lemma 4.4] there exists a level surface $P_\lambda$
which is neither $\r L$ nor $\r H$ with respect to $Q_\mu$.
The pair of surfaces
$(P_\lambda,
Q_\mu)$ gives then rise to a pair of intersection graphs $(G,H)$
without trivial loop.
Since $S^3$ does not contain a $3$-sub-manifold with non trivial torsion,
$H$ does not represent all types.
Therefore $G$ contains a \cs, by Theorem~3.1, which implies
that $Q_\mu$ is $\r C$, for all $\mu<\mu_0$.
\cqfd

{\tmg L{emma} 4.3}\qua (Extremal conditions)
\qua
{\sl Without loss of generality, we may suppose that:
\items
\item[\rm(i)] $P_{\lambda=0}$ is $\cal H$ and $P_{\lambda=1}$
is $\cal L$; and

\item[\rm(ii)]$Q_{\mu=0}$ is $\cal C$ and $Q_{\mu=1}$ is $\cal L$.
\enditems}

\proof
By the previous lemma, the level surfaces $Q_\mu$ in a
thin neighbourhood of $\Sigma$ must be $\cal C$ so is $Q_{\mu=0}$.
Modulo isotopy, we can say
that the level surfaces $Q_\mu$ in a neighbourhood of $Q_{\mu'}$ are
$\cal L$ so is $Q_{\mu=1}$. In the same way of isotopies, we have
the conditions (i).
\cqfd

{\tmg L{emma} 4.4}
\items\sl
\item[\rm(i)] If $K_\gamma$ is a non-trivial knot then for all $\lambda$,
there exists $\r
X_\lambda\in\{\r H, \r L\}$ such that, for all $\mu$, $P_\lambda$
is either $\r N$ or $\r X_\lambda$ with
respect to $Q_\mu$.\par
\item[\rm(ii)] If $s\not=1$ then for all $\mu$,
there exists $\r X_\mu\in\{\r L, \r C\}$ such that, for all $\lambda$,
$Q_\mu$ is either $\r N$ or
 $\r X_\mu$ with respect to $P_\lambda$.
\enditems

\proof (i)\qua
If there exists $(\lambda,\mu)$ such that $P_\lambda$ is $\r C$
with respect to
$Q_\mu$ then the graph $G\subset\ch Q_\mu$ of the associated pair of
intersection
graphs $(G,H)$, contains a \cs, which implies that $S^3$ contains a
(non-trivial)
lens space as a connected summand [8, Lemma~2.5.2.b].
Therefore for all $(\lambda,\mu), \ P_\lambda$ is never $\r C$ with
respect to  $Q_\mu$.

Now, assume that there exist $\lambda,\mu_1$ and $\mu_2$ such that
$P_\lambda$ is $\cal H$ with respect to $Q_{\mu_1}$
and $\cal L$ with respect to $Q_{\mu_2}$.
Then, by [22, Lemma 1.1]
if $K_\gamma$ is not trivial, the low and high disks give an isotopy on
$K_\gamma$ that leads to a minimization of the complexity, which is a
contradiction.

(ii)\qua
Since $\mu<\mu',\ Q_\mu$ cannot be $\cal H$.
So, assume there exist $\lambda_1,\lambda_2$ and $\mu$ such that
$Q_\mu$ is $\cal C$ with respect to $P_{\lambda_1}$ and
$\cal L$ with respect to $P_{\lambda_2}$.
Let $b$ denote the arc component of the intersection of the high disk
(in $P_{\lambda_2}$) with $Q_\mu$ (see
Figure~4(a)) and $\{e_i\}_{1\leq i\leq n}$ the arc components of
the intersection between
the carrying disk (in $P_{\lambda_1}$) and $Q_\mu$ (see Figure~5(a)).
Since there is no local minimum of level lower than $\mu'$,
the arc $b$ can be isotoped in $X$, to $\Sigma$.
The isotopy define a disk $D_{\mu'}$ such that $\mu'$ corresponds to
an inessential local maximum because $b$ does not intersect the carrying
disk and so $D_{\mu'}$ cannot be a singular disk.
But this is a contradiction.
\cqfd

\sh{Graph of singularities}

From now on, we suppose that $K_\gamma$ is not trivial, and
$s\not=1$. For convenience, we say that a  level surface $P_\lambda$ or
$Q_\mu$
is at most $\cal H, \cal L$ or $\cal C$, in reference to the previous lemma.
Note that $P_\lambda$ is never $\r C$ and $Q_\mu$ is never $\r H$.

Let $\Gamma$ be the {\it Cerf graph of singularities},
representing the singularities which appear
in the intersection of two one-parameter families of
surfaces (Figure~7).

A point in $\Gamma$ is a couple of parameters
$(\lambda,\mu)\in[0,1]^2$ for
which the corresponding surfaces $P_\lambda$ and $Q_\mu$ are tangents;
a point
in the exterior of the graph $\Gamma^c=[0,1]^2-\Gamma$, corresponds to
transverse surfaces.

\stdfigure 80mm by 314mm (07 scaled 200) (Figure~7:
Graph of singularities)

Only two types of singular points can appear in $\Gamma$:

\begin{description}
\item[Index~$1$ point]which corresponds to interchange of two
tangency points.
For example, the surface $P_{\lambda=i-\varepsilon}$ in Figure~7, has two
tangency points on two different levels (for two surfaces $Q_{\mu_1}$ and
$Q_{\mu_2}$). Now, decreasing $\varepsilon$ through $0$,
these two points pass
through the same level as tangency points between $P_{\lambda=i}$ and some
surface $Q_\mu$.

\item[Index~$2$ point]which corresponds to birth/death of two
tangency points.
For example, the surface $P_{\lambda=j-\varepsilon}$ in Figure~7,
has two
tangency points on two different levels which degenerate
in a single tangency
point for $\varepsilon$ decreasing to $0$ and disappears for $\lambda>j$.
\end{description}

 All such singular points of $\Gamma$ can be supposed,
 without loss of generality, with distinct abscissa $\lambda$,
and also with distinct
ordinates $\mu$.

Moreover, the slope of such a graph can be supposed neither vertical nor
horizontal. These conditions can be realized by transversality
arguments due to Cerf~[7, Chapter~2].
Furthermore, by isotopies on $\{P_\lambda\}_{\lambda\in[0,1]}$,
the extremal conditions (Lemma~4.3) continue to hold.
Note that all conditions we work with, are open ones.

\pv
{\tmg L{emma} 4.5}
\qua
{\sl
For all $(\lambda,\mu)$ in a connected component of $\Gamma^c$, all the
$P_\lambda$ have the same characteristic in $\{\r H,\r L,\r N\}$ with
respect to $Q_\mu$; and similarly,
all the $Q_\mu$ have the same characteristic in $\{\r C,\r L,\r N\}$ with
respect to $P_\lambda$.
}

\proof
Because $P_\lambda\cap Q_\mu$ cannot change its isotopy
class, except as passing through a point of the graph $\Gamma$.
\cqfd

So, we can associate to each component of $\Gamma^c$
two characteristics from
the set $\{\cal H, \cal L, \cal C, \cal N\}$:
one with respect to $\lambda$
and the other, with respect to $\mu$.

From Lemmas~4.4 and~4.5, we then obtain:

\pv
{\tmg C{orollary} 4.6}
\qua
\items
\sl
\item[\rm(i)] Let $\lambda_1\in [0,1]$ and $\Lambda_1$ the set of the connected
components of $\Gamma^c$ intersecting the vertical line
$\lambda=\lambda_1$.
Then there exists $\r X\in\{\r H, \r L\}$ such that, for all
$(\lambda,\mu)$ in $\Lambda_1$,
$P_\lambda$ have the  characteristic $\r X$ or $\r N$ with respect to $Q_\mu$.

\item[\rm(ii)] Let $\mu_1\in [0,1]$,
and $\Lambda_1$ the set of the connected
components of $\Gamma^c$ intersecting the horizontal line
$\mu=\mu_1$.
Then there exists $\r X\in\{\r L, \r C\}$ such that, for all
$(\lambda,\mu)$ in $\Lambda_1$,
$Q_\mu$ have the  characteristic $\r X$ or $\r N$ with respect to $P_\lambda$.
\enditems

\eject
{\tmg L{emma} 4.7}

\items\sl
\item[\rm(i)] $\forall\lambda\in[0,1]\ \ \exists\mu\in[0,1]\ \vert\
(\lambda,\mu)\in\Gamma^c$ and $Q_\mu$ is $\cal N$ with respect to
$P_\lambda$.\par
\item[\rm(ii)] $\forall\mu\in[0,1]\ \ \exists\lambda\in[0,1]\ \vert\
(\lambda,\mu)\in\Gamma^c$ and $P_\lambda$ is $\cal N$ with respect
to $Q_\mu$.
\enditems

\proof (i)\qua
Assume for a contradiction that there exists $P_{\lambda_1}$,
such that $Q_\mu$ is either $\cal L$ or
$\cal C$ with respect to $P_{\lambda_1}$,
for all $\mu\in [0,1]$ with $(\lambda_1,\mu)\in\Gamma^c$.
By Lemma 4.4, $Q_\mu$ cannot be $\cal L$ and $\cal C$.
Therefore,
there exists a saddle tangency level, $\mu_1$, in $P_{\lambda_1}$,
which is above  the Carrier spheres $Q_{\mu<\mu_1}$
and below the Low spheres $Q_{\mu>\mu_1}$.

Let $\varepsilon>0$ and $\dps
\bigcup_{i=\lambda_1-\varepsilon}^{i=\lambda_1+\varepsilon}P_i$
a regular neighbourhood of $P_{\lambda_1}$.
Since $\mu_1$ corresponds  to a saddle point,
$Q_{\mu_1}$ is $\cal L$ and $\cal C$ with respect to
$P_{\lambda_1-\varepsilon}$ and $P_{\lambda_1+\varepsilon}$, respectively;
in contradiction with Lemma~4.4.

(ii)\qua
We apply the same argument, changing the   Carrier characteristic
to the High one.
\cqfd

{\tmg C{laim} 4.8}
\qua
{\sl Without loss of generality,
we may suppose that there is no pair $(\lambda,\mu)\in\Gamma^c$
such that $P_\lambda$ and $Q_\mu$ are both $\r N$ one with respect
to the other.}

\proof
Otherwise the associated pair of intersection graphs is in contradiction
with the Theorem 3.1.
\cqfd

Therefore, in a same connected component of $\Gamma^c$, we cannot have both
characteristics of $P_\lambda$ and $Q_\mu$ to be $\cal N$.

Let $t=\sup\{\mu\in[0,1]\mid Q_\mu\ {\tm {is}}\ \cal C\}$;
so $0<t<1$, by Lemma~4.3.
Since the slopes of $\Gamma$ are non-horizontal,
the corresponding point $(s,t)\in\Gamma$ is a singular point.

If $(s,t)$ is an index~$2$ point then there are exactly two
connected components of $\Gamma^c$ in a neighbourhood of $(s,t)$.
Therefore, from Corollary~4.6(ii), the level surfaces $Q_\mu$
have characteristic $\cal C$ in one of them and $\cal N$ in the other,
with respect to all the corresponding $P_\lambda$'s.

Furthermore,
$Q_{\mu=t+\varepsilon}$ is $\cal N$ with respect to $P_\lambda$ for all
$\lambda\in [0,1]$; otherwise, as $t=\sup\{\mu\in[0,1]\mid Q_\mu\ {\tm
{is}}\ \cal C\}$ we could find $\lambda_1, \lambda_2$ such that
$Q_{\mu=t-\varepsilon}$ is $\r C$ for $P_{\lambda_1}$
and $Q_{\mu=t-\varepsilon}$ is $\r L$ for $P_{\lambda_2}$,
in contradiction with Lemma~4.4(i).
By  Lemma~4.7, we can find $P_\lambda$ and
$Q_\mu$, each being $\cal N$ with respect to the other;
which is impossible by the previous claim.

Consequently, we may assume that $(s,t)$ is an index~$1$ point.

\stdfigure 60mm by 199mm (08 scaled 250) (Figure~8:
Configuration of {$(s,t)$} in $\Gamma$)

\pv
{\tmg L{emma} 4.9}
\qua
{\sl
The graph $\Gamma$,
in a neighbourhood of the point $(s,t)$ is described in Figure~8 with the
following properties:
\items

\item{} $(\lambda,\mu)\in R_1\Rightarrow P_\lambda$ is $\cal H$
with respect to
$Q_\mu$ and $Q_\mu$ is $\cal N$ with respect to $P_\lambda$.

\item{} $(\lambda,\mu)\in R_2\Rightarrow P_\lambda$ is $\cal N$
with respect to
$Q_\mu$ and $Q_\mu$ is $\cal C$ with respect to $P_\lambda$.

\item{} $(\lambda,\mu)\in R_3\Rightarrow P_\lambda$ is $\cal L$
with respect to
$Q_\mu$ and $Q_\mu$ is $\cal N$ with respect to $P_\lambda$.

\item{} $(\lambda,\mu)\in R_4\Rightarrow P_\lambda$ is $\cal N$
with respect to
$Q_\mu$ and $Q_\mu$ is $\cal L$ with respect to $P_\lambda$.
\enditems
}

\proof
The definition of $t$ implies
$Q_\mu$  is $\r C$
with respect to $P_\lambda$, for $(\lambda,\mu)\in R_2$ and
$Q_\mu$  is not $\r C$
with respect to $P_\lambda$, for $(\lambda,\mu)\in R_4$;
so is $\r L$ or $\r N$ in $R_4$.

As singular points of $\Gamma$ are on different ordinates,
$(s,t)$ is the only singular point of
$\Gamma$ in $N(\mu=t)=[0,1]\times [t-\varepsilon,t+\varepsilon]$
for small enough $\varepsilon$.
So, as $Q_\mu$ is $\r C$ in $R_2$, Corollary~4.6(ii) implies:
\smallskip

\centerline{$Q_\mu$ is $\r N$ with respect to $P_\lambda$ for all
$(\lambda,\mu)\in N(\mu=t)\backslash(R_2\cup R_4)$\ \ $(*)$}

Furthermore, if we suppose $Q_\mu$ is $\r N$ in $R_4$ then, by Lemma~4.7(ii),
we deduce a contradiction to Claim~4.8.
Then $Q_\mu$ is $\r L$ in $R_4$ and $\r N$ in $R_1\cup R_3$.

So, Claim~4.8 again implies $P_\lambda$ is not $\r N$ in $R_1\cup R_3$ and
$(*)$ implies $P_\lambda$ is not
$\r N$ in $N(\mu=t)\backslash (R_2\cup R_4)$.
And from Lemma~4.3 (Extremal conditions),
we conclude $P_\lambda$ is $\r H$ and $\r L$ in $R_1$ and
$R_3$, respectively.
Finally, Corollary~4.6(i) implies $P_\lambda$ is $\r N$ in $R_2\cup R_4$.
\cqfd

For convenience, we note $\lambda=s\pm\varepsilon$, abscissa of points in
 $\Gamma$, in a regular neighbourhood of $s$. And similarly for
ordinates $\mu=t\pm\varepsilon$, in a neighbourhood of $\mu=t$.
Note that (see Figure~8) there exist $\varepsilon'>0$ and $C_\lambda$,
$C'_\lambda$ with $t-\varepsilon'<C_\lambda<t<C'_\lambda<t+\varepsilon'$
such that:

\items

\item[\rm(i)] $P_{\lambda=s-\varepsilon}$ and $Q_{\mu=C_\lambda}$
are tangents in a point $x_\lambda$.

\item[\rm(ii)] $P_{\lambda=s-\varepsilon}$ and $Q_{\mu=C'_\lambda}$ are tangents
in
a point
$x'_\lambda$.

\item[\rm(iii)]$x_\lambda$ and $x'_\lambda$ are in the boundary of $R_1$
(on different lines of $\Gamma$).
\enditems

Therefore, the point $(s,t)\in\Gamma$ is an index~$1$ point, which
corresponds to a pair of saddle-points $x=x_{\lambda=s}$ and
$x'=x'_{\lambda=s}$
which are the two tangency points between the surfaces $P_{\lambda=s}$ and
$Q_{\mu=t}$.
For convenience, the reader should refer, say [22, p.381].

Let $D^-$ and $D^+$ be the limits in $Q_{\mu=t}$ when $\lambda$ goes to
$s$, of the  low and high disks for $P_\lambda$ coming from $R_1$ and
$R_3$, respectively.
Then the boundary of $D^-$ (resp.\ $D^+$)
contains $x$ and $x'$.

Indeed, if for example $x\notin\partial D^-$ then the
low disk for $(\lambda,\mu)\in R_1$ does not disappear when $(\lambda,\mu)$
goes through the line defined by $\{\mu=C_\lambda\}$ (Figure~8). Therefore
$P_\lambda$ is still high with respect to $Q_\mu$ for $(\lambda,\mu)\in R_2$
and this is in contradiction with Lemma 4.7.

If $x'\notin\partial D^-$ then we pass the low disk of $R_1$ through the
line
defined by $\{\mu=C'_\lambda\}$ (Figure~8),
surviving so in $R_4$ in contradiction with Lemma 4.7.

For the case $x\notin\partial D^+$, one can pass the high disk from $R_3$ to
$R_4$ and for the case $x'\notin\partial D^+$, from $R_3$ to $R_2$; arriving
also at a contradiction.

Trying to put the limit-disks $D^-$ and $D^+$ in $Q_{\mu=t}$, we deduce the
configuration on $Q_{\mu=t}$ of Figure~9(a); otherwise
$P_{\lambda=s}$ has only two boundary components, and hence $K_\gamma$ is
trivial in $S^3$.
Let $V_1$ and $V_2$ denote the boundary components of $Q_{\mu=t}$ that
$D^+$ and $D^-$, respectively meet (Figure~9(a)); see [22, pp. 383-384],
for  convenience.

\begin{figure}[ht!]
\cl{\includegraphics[width=80mm]{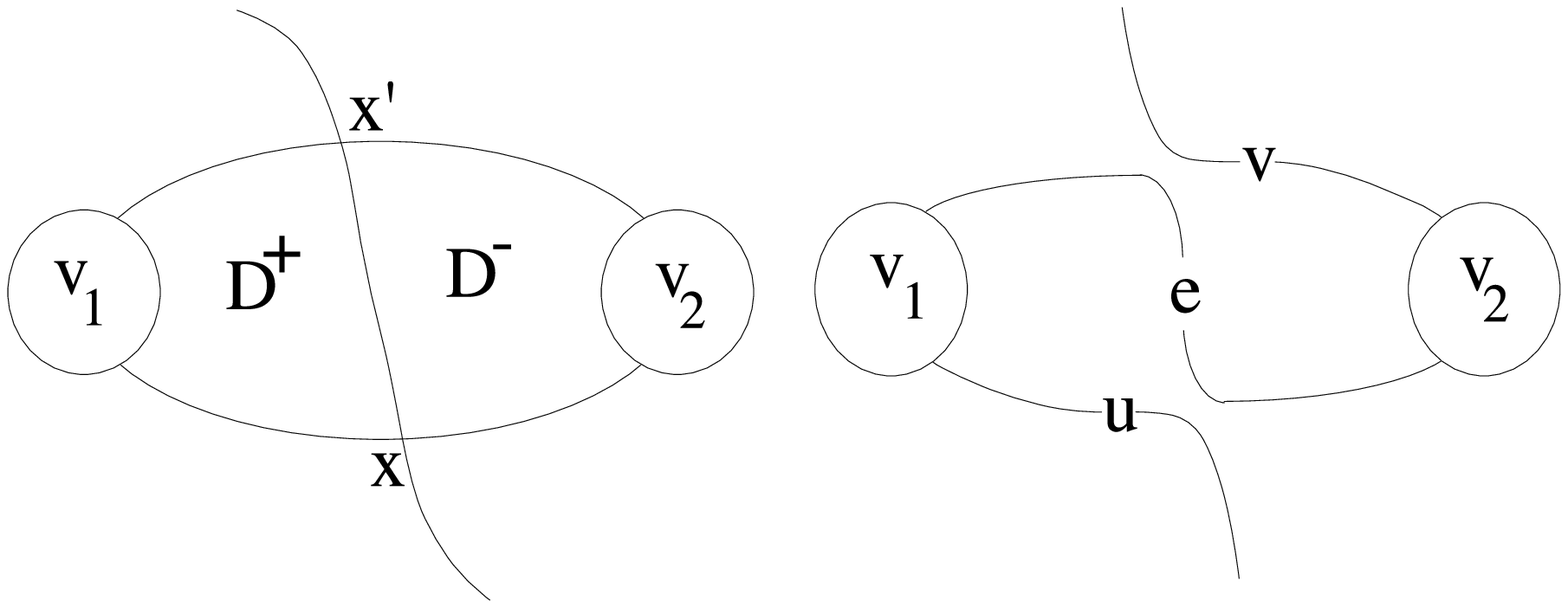}}
\medskip
\small\qquad\qquad\qquad\qquad (a)\qua $D^\pm$ in $Q_{\mu=t}$ \qquad\qquad
(b)\qua $b_\varepsilon$
in
$Q_{\mu=t+\varepsilon}\cap P_{\lambda=s}$
\nocolon\caption{}
\end{figure}

Now we look at the characteristics $\r C$ and $\r L$ of $Q_\mu$,
with respect to $P_\lambda$, that is for $(\lambda,\mu)$ in $R_2$ and $R_4$,
respectively.
The characteristic $\cal L$ in $R_4$ gives rise to a high disk
$D_\varepsilon$ ($\subset P_{\lambda=s}$)
whose boundary contains an arc $b_\varepsilon$ in $Q_{\mu=t+\varepsilon}$,
i.e.\
$b_\varepsilon=\partial D_\varepsilon\cap Q_{\mu=t+\varepsilon}=
D_\varepsilon\cap Q_{\mu=t+\varepsilon}$.
Let $b^+$ be its limit for $\varepsilon$ decreasing to $0$.
Then $b^+$ is a simple arc properly embedded in $Q_{\mu=t}$,
which contains $x$ and $x'$ for the same reason as above with
$\partial D^-$ (and $\partial D^+$).

Furthermore, $x$ and $x'$ are saddle-points, so we describe in Figure~9(b)
the
local intersection in $Q_{\mu=t+\varepsilon}$ with $P_{\lambda=s}$.

Let $e$ be the  arc in $Q_{\mu=t+\varepsilon}$, joining $V_1$ and $V_2$,
which
corresponds to the arc in $Q_{\mu=t}$ joining
$V_1$ and $V_2$, and passing through $x'$ and $x$ (Figure~9(a)).
Let $u$ and $v$ be the remaining open arcs.
Then either $b_\varepsilon$ contains $e$ or $u\cup v$.
Consequently, $b^+$ joins the two components $V_1$ and $V_2$.

Likewise, the characteristic $\cal C$ in $R_2$ gives an union of arcs in
$Q_{\mu=t-\varepsilon}$. Their limit $b^-$ (union of arcs)
must also contain $x$ and $x'$.
Furthermore, $V_1$ and $V_2$ are joined by $b^-$,
as well as for $b^+$ above.

Thus the \nd\ $K$ intersects exactly twice the sphere $Q_{\mu=t}$
and once the spine $\Sigma$. So $s=1$, proving Theorem~4.1.

\section{Geometric intersection with a lens spine}

Let $K$ be a knot in  a lens space $L$,
which produces $S^3$ by $\gamma$-Dehn surgery.
Let $\Sigma$ be a standard spine of $L$, and
$s$ be the minimal geometric intersection number between $K$ and $\Sigma$.
Let $K_\gamma$ be the core of the surgery in $S^3$.

A {\it $(m,n)$-cable} of a knot $k$, is a knot
lying on the boundary of a regular neighbourhood of $k$,
which goes $m$ times in the meridional direction and
$n$ times in the longitudinal.
And non-trivial cables of torus knots are known to produce
lens spaces by Dehn surgery [2, 20].
Note that cables of the trivial knot (or
$0$-bridge braids) in $S^3$ are exactly the torus knots.

\pv
{\tmg L{emma} 5.1}
\qua
{\sl Assume that $s=1$.

\items
\item[\rm(i)] If a thin presentation of $K$,
with respect to $\Sigma$,
begins by a local maximum then $K$ is a $0$
or $1$-bridge braid in $L$.
\par
\item[\rm(ii)] If a thin presentation of $K$,
with respect to $\Sigma$,
begins by a local minimum then $K_\gamma$ is a cable in $S^3$.
\enditems}

\proof
(i)\qua
The knot $K$ intersects exactly twice the sphere $S=\partial N(\Sigma)$.
Let $X$ be the $3$-ball bounded by $S$, and $a=K\cap X$, then $(a,X)$ is a
trivial tangle, since $a$ is a ``local maximum''.
Therefore there is a disk, $D$ say, in $X$
such that $\partial D=a\cup b$, where $b$ is a simple arc in $S$.
Then, we can isotope $a$, successively via $D$ and $N(\Sigma)$,
to a simple arc
$\alpha$ in $\Sigma$, with its endpoints identified, since $s=1$.

 Let $\kappa$ be the core of the spine.
Then $\alpha$ and
$\kappa$ ``cobound'' a pinched spinal helix of order less than
$\hbox{\ssmall order of }\Sigma\over2$ (Figure~10).
In other words, $K$ is isotoped to a simple closed curve in $\Sigma$ which
intersects the core $\kappa$ only once.
Thus, $K$ is isotoped to $\alpha$ and therefore, is a $0$
or $1$-bridge braid in $L$.

\begin{figure}[ht!]
\cl{\includegraphics[width=85mm]{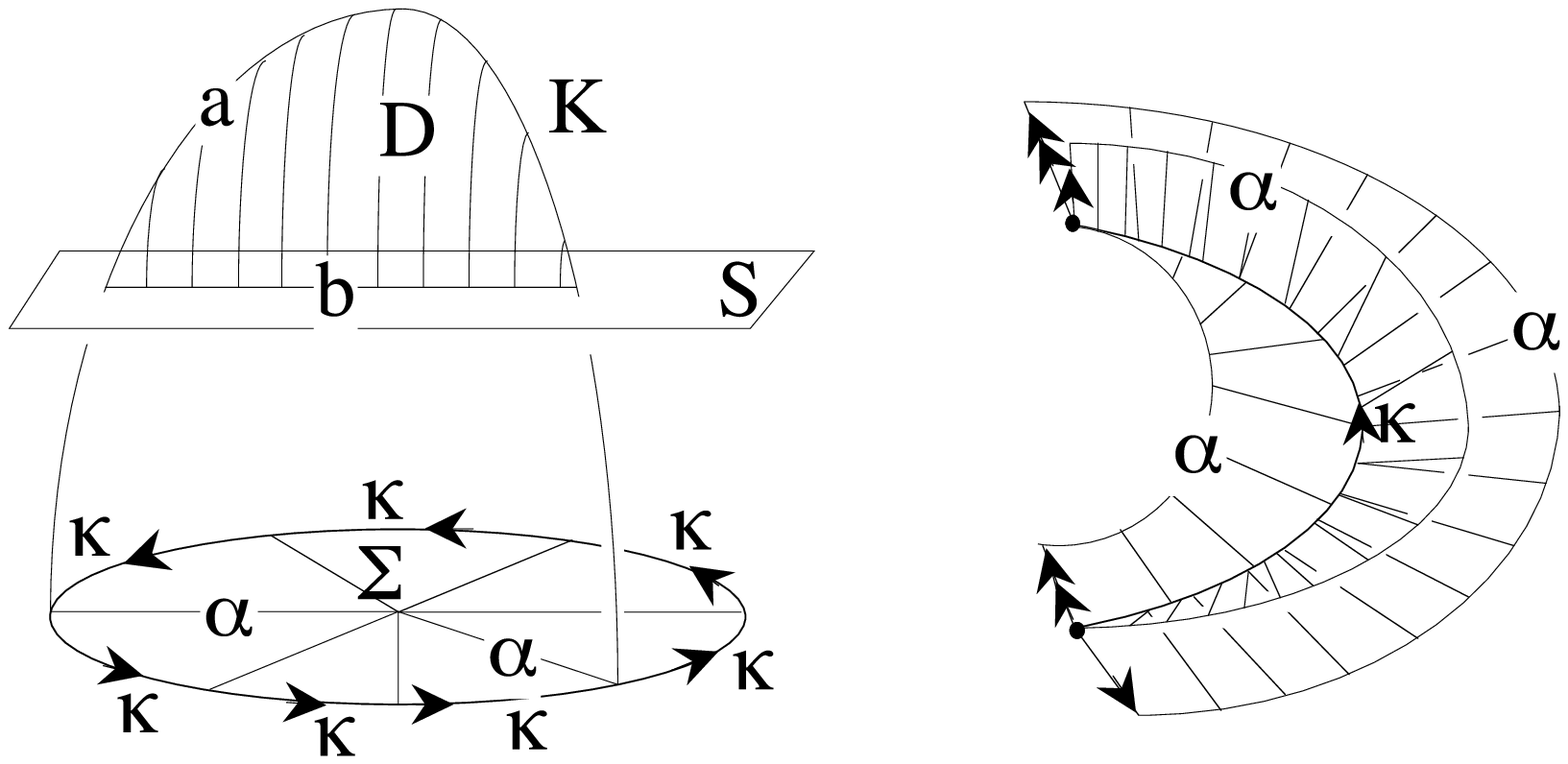}}
\medskip
\small\qquad\qquad\qquad (a)\qua High disk until $\alpha\subset\Sigma$
\qquad\qquad(b)\qua Pinched spinal
 helix 
\nocolon\caption{}
\end{figure}

(ii)\qua
Let $\widehat Q=\partial N(\Sigma)$ be
a level $2$-sphere of the thin presentation of $K$,
in a neighbourhood of the spine.
And denote $Q=\widehat Q\cap L_K$, the corresponding level surface.
Since $\vert\widehat Q\cap K\vert=2$ and
the thin presentation of $K$ begins by a local minimum,
we deduce that $Q$ is incompressible and
$\partial$-incompressible in $L_K$.
Fixing a thin presentation of $K_\gamma$ in $S^3$,
$Q$ is a properly embedded planar surface in $X_{K_\gamma}$ such that
$\partial Q$ is not meridional.
Therefore, by [15, Lemma~4.4],
one can find a level surface $P$ from the $S^3$-foliation,
transverse to $Q$ in $L_K$ with the property that
$P$ (and already $Q$) is neither $\r L$ nor $\r H$
with respect to $Q$ (resp.\ $P$).

We can then associate a pair of intersection graphs $(G,H)$ of type
$(\ch P,\ch Q)$ as in Section~2 such that
$G$ and $H$ contain no trivial loop.
Now, $H$ has exactly two vertices, therefore the edges in $H$,
all join these two with distinct labels at their endpoints.
So, $K_\gamma$ is a cable knot by [21, Section~5].
\cqfd

Now Theorem 4.1 and Lemma 5.1(i) prove together Theorem 1.2.
Indeed, if we have the hypothesis of Theorem~1.2,
so a knot $K$ in an essential thin
presentation beginning by a
local maximum then Theorem~4.1 applies and Lemma~5.1(i) implies the result.

\sh{Proof of Proposition 1.3}

Let us suppose, in this paragraph
that $K$ is neither a $0$ nor $1$-bridge braid in $L$.
And furthermore, let us suppose that $s<3$.

First, if $s=1$ then $K_\gamma$ is a cable knot in $S^3$ by
Lemma~5.1(ii).
So, there is a
non-meridional essential annulus $A$
in the exterior $X_{K_\gamma}=X_K$.
We also have the level surface $Q$, below the first local minimum,
 which is
an essential annulus in $X_K$ transverse to $A$.
We obtain a pair of graphs, each with two vertices and
without trivial loop.
And this contradicts the parity rule.

Now, if $s=2$ and the thin presentation of $K$ is beginning
by a local maximum then it is inessential, by Theorem~1.2.

\pv
{\tmg C{laim} 5.3}
\qua
{\sl If $K$ is standardly spinal then it is a $0$ or $1$-bridge
braid.}

\proof
As $s=2$, then $K=\alpha\cup\beta\subset\Sigma$ such that
$\alpha\cap\kappa=\beta\cap\kappa=\partial\alpha=\partial\beta$,
where $\kappa$ is the core of $\Sigma$.
So, $K$ is a $0$ or $1$-bridge braid in $L$.
\cqfd

This claim implies that the thin presentation of $K$ has a local
minimum above the first and inessential local maximum.
So, by the same argument used in the proof of Lemma~5.1(ii), we deduce
that the core of the surgery is a cable knot in $S^3$.
And using the same reasonning as above for the case $s=1$, we arrive
at a contradiction with parity rule.

So, from now on, we may suppose that for $s=2$, we only have the case
in which the thin presentation of $K$ begins by a local minimum.
Then let \ch Q be the level $2$-sphere
in a neighbourhood of
the spine $\Sigma$; i.e.\ such that
the corresponding level surface $Q=\partial(N(\Sigma)\cap L_K)$.

If $Q$ is not incompressible, then there is a compressing disk
in $L_K$ that
 cobounds a $3$-ball with a disk $D$ in \ch Q;
$IntD$ intersects then exactly twice the knot $K$,
for its boundary to be
essential in $Q$.
We isotope $K$, using this $3$-ball, below the first local minimum,
reducing
so the complexity of the thin presentation of $K$,
which is supposed to be
impossible.
So, $Q$ is incompressible in $L_K$.
Furthermore, because of the minimal complexity of this
presentation beginning
by a local minimum, $Q$ is also $\delta$-incompressible.

Then, by [15, Lemma~4.4], we find a level sphere \ch P in
a thin presentation of $K_\gamma$ in $S^3$ such that
$P$ is incompressible and $\partial$-incompressible in $X_{K_\gamma}=L_K$.

So, the pair of graphs $(G,H)$ of type $(\ch P,\ch Q)$ does contain
no trivial
loop.
Applying Theorem~3.1, $G$ then contains a \cs\ $\sigma$
of length $n$, say;
let
$\{1,2\}$ be the labels of $\sigma$ and $W_1$, $W_2$
the corresponding vertices
in $H$.
If $\vert G^0\vert=p$ and $\Delta=\Delta(\infty,\gamma)$, then there are
$\Delta p$ endpoints of edges around each vertex of $H$.
The edges of $\sigma$ separate $H$ in $n$ disjoint bigons if we consider
the vertices as points;
i.e.\ $\widehat Q-N(\sigma\cup W_1\cup W_2)$ is a union of $n$ disks.

\stdfigure 67mm by 152mm (11 scaled 430) (Figure~11:
Standard spine in $Z$)

\pv
{\tmg C{laim} 5.4}
\qua
{\sl The two other vertices of $H$ ($W_3$ and $W_4$ say) are in
a same bigon $B$.}

\proof
Suppose they do not. Since $H$ contains no trivial loop,
all endpoints of $W_i$
($i\in\{3,4\}$) are incident to edges meeting only $W_1$ or $W_2$.
Therefore, $W_1$ and $W_2$ are incident to more than $2\Delta p$ edges
(counting the edges of $\sigma$), which is impossible.
\cqfd

 So $B^*=\ch Q-$int$B$ contains  $W_1$, $W_2$ and
the edges of $\sigma$.
Let $J_{1,2}$ be the $1$-handle in $N(K)$, bounded by $W_1$ and $W_2$,
which does not contain $W_3$ (and therefore $W_4$).
We consider a regular neighbourhood $W$ of $B^*\cup J_{1,2}$ which is
a solid torus whose boundary is pierced twice by $K$.
We then add the disk-face $D$ bounded by
$\sigma$, as a $2$-handle: $Z=N(W\cup D)$ is a punctured lens space of
order $n$.
Note that $D$ is not intersected by $K$ because $D$ is a disk-face.
But $\partial D$ (which is $\sigma$) is the boundary of a spinal helix
of order $n$ in the punctured lens space $Z$:
the core $\kappa$ of the solid torus $W$ is the singular set (the core)
of a $2$-complex bounded by $\partial D$.

\stdfigure 72mm by 201mm (12 scaled 400) (Figure~12:
Isotopy of $\tau$
with respect to the spinal helix)

Let $\kappa_1$ and $\kappa_2$ be two arcs such that
$\kappa=\kappa_1\cup \kappa_2$ and $\kappa_1$ is isotoped to
$K\cap W=\tau$
(see Figure~11).
By capping off this spinal helix with $D$,
we then obtain a standard spine
$\Sigma$ of $Z$.
Then we may isotope $\tau$, fixing $\partial\tau$
(see Figure~12) to intersect
$\Sigma$ in a single point in $\kappa$.
So $K$ intersects a standard spine exactly once, which is a contradiction
with Lemma~5.1.

\section{Proof of Theorem 1.4}

Let $K$ be a non-standardly spinal knot in a lens space
$L=L(m,n)$ with standard spine $\Sigma$,
which yields $S^3$ by $\gamma$-Dehn surgery.
A fortiori, $K$ is not a $0$ or $1$-bridge braid in $L$.
Denote $s$ the minimal
geometric intersection number between $K$ and $\Sigma$.
Let $k=K_\gamma$ be the core of the $\gamma$-Dehn surgery in $S^3$.
As we remark after Proposition~1.3 in the Introduction,
we may assume that $k$ is not a torus knot.
So there is
an integer slope $r$ by~[8] such that $r$-Dehn surgery on $k$ yields $L$.
Then $H_1(L)=\d Z/\vert r\vert\d Z$ and
$\Delta=\Delta(r,\lambda)=\vert r\vert=m\geq 2$.

Furthermore, since $K$ is neither $0$ nor $1$-bridge braid,
we have $s\geq3$ (Proposition~1.3), and hence by Theorem~4.1,
we may assume that the essential thin presentation of $K$ begins
by a local minimum.

So, let \ch Q be a level $2$-sphere between $\Sigma$ and
the first local minimum and note $Q=\widehat Q\cap L_K$
the corresponding level surface.

And let $P$ be a Seifert surface (with minimal genus) for $k$,
properly embedded
in $L_K$.
The single boundary component of $P$ has the longitudinal slope
$\lambda$ and
let denote \ch P, the closed surface obtained by
capping off $P$ by a disk.

Now, we consider the pair of intersection graphs $(G,H)$ of type
$(\ch P,\ch Q)$.

\pv
{\tmg C{laim} 6.1}
\qua
{\sl The graphs  $G$ and $H$  contain no trivial loop.
}
\proof
Since $P$ is  an incompressible and $\partial$-incompressible
surface, then $H$ contains no trivial loop.
And if $G$ contains a trivial loop, then we may assume (up to isotopy)
that the
tight presentation of $K$ begins by a local maximum,
which is a contradiction.
\cqfd

Let $G_x$ be the subgraph of $G$, with the single vertex and
all the edges with
an endpoint labelled by $x$.
By the parity rule, the edges of $G_x$ do not have both endpoints
labelled by $x$.
\pv
{\tmg C{laim} 6.2}

{\sl If $\Delta\geq 2g+1$ then $G_x$ contains at
least two disk-faces.}

\proof
The Euler characteristic calculation for $G_x$ gives $\chi(\ch
P)=2-2g=V-E+F$, where $V$ is the number of vertices, $E$ is the number
of edges of $G_x$, and $\dps F=\sum_{f\ face \ of G_x}\chi(f)$.
Here we have $V=1$ and $E=\Delta$, so
$F=\Delta-2g+1$. Therefore, if $\Delta\geq 2g+1$ then $G_x$ contains at
least
two disk-faces.
\cqfd

Each disk-face of $G_x$ contains a \cs\ in $G$, by [23, Lemma~2.2].
If $G$ contains a \cs\ of length $\ell$, then
$X_k(r)$ contains a lens space of order $\ell$ (see [8] or below);
therefore $\ell=m=\vert r\vert=\Delta$.

Assume for a contradiction that $\Delta\geq 2g+1$.
Let $\sigma$ and $\sigma'$ be two \cs s in $G$.
Without loss of generality, we may assume that $\{1,2\}$ are the labels of
$\sigma$.
Since $\sigma$ contains $\Delta$ edges, the corresponding edges in $H$ join
the vertices $W_1$ and $W_2$, forming a connected component of $H$;
i.e.\ no more edges than the $\Delta$'s from $\sigma$, can have its
extremities attached on $W_1$ or $W_2$.
Therefore, $\sigma'$ has a pair of labels $\{n,n+1\}$ disjoint from
$\{1,2\}$.
And the vertices $W_n$, $W_{n+1}$, together with the corresponding edges of
$\sigma'$ form another connected component of $H$.

Let $\Lambda$ and $\Lambda'$ be these components of $H$, corresponding to
$\sigma$ and $\sigma'$ respectively.
Since $H$ is on a $2$-sphere, there exist two disjoint disks $D$ and $D'$ in
\ch Q, containing $\Lambda$ and $\Lambda'$ respectively.

Therefore $L=X_k(r)$ contains two disjoint punctured lens spaces.

Indeed, let $A, A'$ be the face disks bounded by
respectively $\sigma, \sigma'$ in $G$;
$V$ the $2$-handle of the attached solid torus between respectively
$W_1$ and $W_2$ with no other $W_i$ inside,
and similarly
$V'$ the $2$-handle of the attached solid torus between respectively
$W_n$ and $W_{n+1}$ with no other $W_i$ inside.
Then $N(D\cup A\cup V)$ and $N(D'\cup A'\cup V')$ are two disjoint punctured
lens
spaces, in contradiction with $L$.

Consequently, $G$ contains at most one \cs.
So $\Delta\leq 2g$.

\Addresses\recd
\end{document}